\documentclass[11pt]{amsart}
\usepackage{amsthm}
\usepackage{amsmath}
\usepackage{amssymb}
\usepackage{hyperref}
\usepackage{todonotes}
\usepackage{tikz}
\usepackage{xcolor}
\usepackage{thmtools}
\usepackage{cleveref}
\usepackage{mathtools}

\definecolor{mypink1}{RGB}{255, 100, 251}

\usetikzlibrary{decorations.pathreplacing}
\usetikzlibrary{arrows.meta, positioning, shapes.geometric}

\usepackage[margin=2.75cm]{geometry}

\newtheorem{theorem}{Theorem}[section]
\newtheorem{proposition}[theorem]{Proposition}
\newtheorem{lemma}[theorem]{Lemma}

\newtheorem*{example}{Example}

\theoremstyle{remark}

\theoremstyle{definition}
\newtheorem{definition}[theorem]{Definition}

\newcommand{\N}{\mathbb{N}}

\newcommand{\R}{\mathbb{R}}

\newcommand{\hess}{\text{Hess}}
\newcommand{\tr}{\text{tr}}
\newcommand{\grad}{\nabla}
\newcommand\inner[2]{\langle #1, #2 \rangle}
\newcommand{\lap}{\Delta}
\newcommand{\Div}{\text{div}}

\newcommand{\ric}{\text{Ric}}
\newcommand{\del}{\partial}

\newcommand{\e}{\varepsilon}
\newcommand{\s}{\Sigma}
\newcommand{\ckaq}{C^{k,\alpha}_{-q}}

\newcommand{\model}{\hat{M}}
\newcommand\moddec[1]{\hat{#1}}

\newcommand{\dist}{\text{dist}}
\newcommand{\p}{\prime}

\newcommand{\bh}{\text{bh}}
\newcommand{\mg}{\hat{g}}

\newcommand{\di}{\delta_1}
\newcommand{\dii}{\delta_2}
\newcommand{\genus}{\mathfrak{g}}

\title{Uniqueness Results for Static Asymptotically Locally Hyperbolic (ALH) Manifolds}
\begin{document}
\author{Katie Gravel}
\maketitle
\begin{abstract}
    We give uniqueness results for related classes of static asymptotically locally hyperbolic (ALH) manifolds. As an application, we give examples of one ended hyperbolic manifolds with boundary that are not optimal for the Positive Mass Theorem or for the Riemannian Penrose Inequality.
\end{abstract}
\section{Introduction} 
Let $n\geq 3$ and let $(M^n,g)$ be a connected orientable Riemannian manifold possibly with boundary. If there is a nonzero smooth function $f\colon M\to\R$ such that 

\begin{equation*}
    L_g^*(f) :=-(\lap f )g + \hess_g f -f\ric(g) = 0 
\end{equation*} we say that $(M,g)$ is \emph{static} and $f$ is a \emph{static potential} on $M$. If the level set $\{f = 0\}$ is nonempty, it is a classical fact that $\{f=0\}$ is a totally geodesic hypersurface. Accordingly, if $M$ has boundary and $f = 0$ on $\del M$, we call the triple $(M,g,f)$ a \emph{static black hole}, where $\del M$ is the `black hole'. For simplicity, throughout the introduction, we assume $\{f = 0\}$ is connected and $f$ has been normalized so that $|\grad f| = 1$ along $\{f= 0\}$.

Since the operator $L_g^*$ is overdetermined, one may expect static manifolds to (1) exhibit interesting geometric properties and (2) be `rare' in some sense. Indeed, given a static manifold $(M,g)$ with static potential $f$, the Lorentzian metric given by \[\bar{g} = -f^2dt^2 + g\] on $\R\times M$ is a solution to the Einstein vacuum equations with cosmological constant $\Lambda = -n(n-1)/2$. This gives a heuristic argument for the strong uniqueness for static metrics: physically, a black hole without rotation or charge should be determined only by its mass. Meaning, we should expect exactly 1-parameter family of asymptotically flat static black holes when $n = 3$. And indeed the only asymptotically flat $3$-dimensional static black holes are the positive mass Schwarzschild metrics (\cite{Is67}, \cite{BuMa87} when the black hole is not assumed to be connected): \[g_{\text{Schw}} = \frac{dr^2}{1-\frac{2m}{r}} + r^2g_{\text{round}},\] on $(2m,\infty)\times S^{2}$ and where $g_{\text{round}}$ is the round metric on $S^2$. For asymptotically flat $n$-dimensional static black holes where $n > 3$, the Schwarzschild-Tangherlini metrics are the unique 1-parameter family of asymptotically flat static black holes (\cite{GiIdSh02}). 

The classification for asymptotically locally hyperbolic (ALH)--see Definition \ref{def:alh}--is more complicated. In this introduction, we focus on the case where $n = 3$. Let $k\in\{0,\pm 1\}$ and let be $(N,b)$ a surface of genus $\genus$ with constant Gaussian curvature $k$. Let $m\in\R$ be large enough so that the largest root $r_m$ of $V(r)^2 = r^2+k-\frac{2m}{r}$ is simple and positive. Then the metric \[g_{\text{Kott}} = \frac{dr^2}{V(r)^2} + r^2b\] on $(r_m,\infty)\times N$ is an ALH static black hole metric, called a \emph{Kottler metric}. For a fixed genus $\genus$, the space of  Kottler metrics with topology $\R\times\s_{\genus}$ is finite-dimensional modulo diffeomorphism because the space of metrics of constant curvature $k$ on $\s_{\genus}$ is finite-dimensional modulo diffeomorphism.

Similar to the asymptotically flat case, partial uniqueness results for Kottler metrics have been established in (\cite{Bo22}, \cite{HaWa26}, \cite{Lee2015}). The situation expected by (\cite{AnChDe02}) when we consider a black hole whose conformal infinity is a hyperbolic surface, $\s_{\genus\geq 2}$, is thought to be more complicated than the asymptotically flat case. For a conformally compact static metric $g$ with static potential $f>0$ such that $f\to \infty$ near infinity, there is a choice of metric at infinity given by \[g_{\infty} :=\lim_{f\to\infty}f^{-2}g.\] It is suggested in \cite{AnChDe02} that there is an infinite dimensional (modulo diffeomorphism) space of metrics $h$ on $\s_{\genus}$ such that there is a static metric $g$ on $\R\times\s_{\genus}$ such that $g_{\infty} = h$. In particular, the space of static metrics with topology $\R\times\s_{\genus}$ is expected to be infinite-dimensional.

We first address the asymptotic behavior of static potentials $f$ on ALH manifolds with conformal infinity $\s_{\genus}$, and more generally, ALH manifolds with conformal infinity of what we call \emph{negative type} (see Definition \ref{def: conf inf}). In particular, the static equations do not imply that $|f|\to\infty$ near infinity in full generality because there are examples of static potentials on $\mathbb{H}^n$ such that $f\to 0$ along some geodesics near infinity. In settings where the space of static black hole metrics is expected to be finite dimensional, we can determine the end behavior of $f$ by inspection. Since we don't expect the space of static ALH metrics with conformal infinity $\s_{\genus\geq 2}$ to be finite dimensional, we must proceed differently. 

We call a geodesic ray $\gamma\colon [0,\infty)\to M$ a \emph{radial geodesic} if $\gamma$ is a geodesic emanating from $\del M_{ext}$ with velocity $\del_r/|\del_r|$ (see discussion after Definition \ref{def:alh} for notation). We prove the following general dimension result:

\begin{restatable}{theorem}{proper}\label{thm: level sur are closed}
    Let $(M,g)$ be a static ALH one-ended black hole, with conformal infinity of negative type and $f> 0$ near infinity. Let $\gamma_x$ be a radial unit speed geodesic emanating from $x\in \del M_{ext}$. Define
\[
F=\{x\in \del M_{ext}:\lim_{t\to\infty}f(\gamma_x(t))=0\}.
\]
Then \(F\) is finite. Moreover, for every nonempty compact set \(K\subset \del M_{ext}\setminus F\),
\[
\lim_{t\to\infty}\inf_{x\in K}f(\gamma_x(t))=\infty.
\]
In particular, \(f\to+\infty\) uniformly on the end if and only if \(F=\emptyset\).
\end{restatable}

Next, since a static ALH black hole $(M,g)$ is a (model for a) black hole in the physical sense, the metric $g_{\bh}$ on the black hole boundary should be sufficient to determine the entire static metric, since the black hole is the only massive body in the system. 

Therefore, one would like to understand the relationship between $g_{\infty}$ and $g_{\bh}$. For the general class of static metrics $g$ on $(0,\infty)\times \s_{\genus}$, there is no established relationship between $g_{\infty}$ and $g_{\bh}$. Therefore, it is a priori possible that two different choices of $g_{\infty}$ could yield distinct static ALH black hole metrics that have the same $g_{\bh}$. To investigate this type of uniqueness, we study the static equations at the black hole. First, we note that by rewriting $L_g^*(f) = 0$ as a coupled system:
\[
    \begin{cases}
        \hess_gf = f(\ric(g)+ng)\\
        \lap_gf=nf,
    \end{cases}
\] 
we see more clearly that the static equations degenerate at $f = 0$. Therefore, analysis at the black hole must be treated carefully in general. Nevertheless, we establish that a choice of metric $g_{\bh}$ cannot determine two static ALH black hole metrics:

\begin{restatable}{theorem}{bonus}\label{thm: bonus}
    Let $n\geq 3$ and let $(N,g_{\bh})$ be a closed $(n-1)$-manifold. Let $M = (0,\infty)\times N$. There is at most one (up to diffeomorphism) ALH static metric $g$ on $M$ such that $g|_{\{0\}\times N}=g_{\bh}$. 
\end{restatable}

Building on Theorem \ref{thm: bonus}, we show a uniqueness result for 2 ended complete ALH $3$ manifolds:

\begin{restatable}{theorem}{main}\label{thm:main}
    Let $(M,g)$ be a complete two ended analytic static ALH $3-$manifold such that $M = \R\times \s_{\genus\geq 2}$. Then, $(M,g)$ is the doubling of an ALH static black hole.
\end{restatable}

As an application of Theorem \ref{thm:main}, we give examples of one ended hyperbolic manifolds with topology $(0,\infty)\times\s_{\mathfrak{g}\geq 2}$, which are not optimal for either the Positive Mass Theorem or the Penrose Inequality. More specifically, the Kottler metrics with hyperbolic conformal infinity and $m = 0$ are exactly hyperbolic. If the mass $m = 0$ Kottler metric is doubled over its boundary, it corresponds to a well-known type of complete hyperbolic manifold called a \emph{Fuchsian} manifold. And, closely related to the Fuchsian manifolds are a family of hyperbolic metrics also on $\R\times\s_{\mathfrak{g}}$ called \emph{quasi-Fuchsian} manifolds (see Section \ref{sect: quasi} for more details on quasi-Fuchsian metrics). Non-Fuchsian quasi-Fuchsian metrics are not static, and therefore, their scalar curvature can be freely perturbed (\cite{corvino2000}, Theorem 1). This free perturbation of scalar curvature leads to the next two theorems of this paper:

\begin{restatable}{theorem}{nopmt} \label{lem: quasi-Fuchsian scal curv}
    Let $(M,g)$ be a quasi-Fuchsian manifold that is not Fuchsian. Let $\Omega\subset M$ be a relatively compact open set, isotopic to a tubular neighborhood of $\{0\}\times \s_{\genus}$. For any sufficiently small, compactly supported $S\colon \Omega\to\R$, there is a metric $g_S$ nearby to $g$ such that $R_{g_S} = -6+S$ and $g \equiv g_S$ outside $\Omega$.
\end{restatable}
Therefore, quasi-Fuchsian manifolds, although hyperbolic, cannot be rigidity cases of a Positive Mass Theorem, since these perturbations should have the same mass as $(M,g)$. 

Second, we address the Penrose inequality. Let $(M^3,g)$ be an ALH manifold with outer-minimizing minimal boundary of area $A$ and $R(g)\geq -6$. Assume that $M$ has conformal infinity of genus $\mathfrak{g}\geq 2$. Let $m(g)$ be the Chrusciel-Herzlich mass of $M$ (see Definition \ref{def: mass}). The Penrose Inequality for this class of manifolds is thought to state: \[m(g)\geq \frac{1}{(\mathfrak{g} -1)^{3/2}}\sqrt{\frac{A}{16\pi}}\Bigg(1-\mathfrak{g} + \frac{A}{4\pi}\Bigg),\] with equality if and only if $(M,g)$ is isometric to a Kottler metric of mass $m(g)$. The Penrose Inequality for ALH $3$-manifolds is proven when the related mass quantity $\bar{m}(g)$ satisfies $\bar{m}(g) \leq 0$ with certain asymptotics on $g$ in \cite{Lee2015}. 

Now, let $(M,g)$ be the one ended hyperbolic manifold obtained by cutting a quasi-Fuchsian manifold along a strictly stable outermost minimal surface, and considering an end that contains no compact minimal surfaces. It can be seen using the Gauss equation that $(M,g)$ satisfies the Penrose inequality strictly. However, we show that $(M,g)$ is not \emph{optimal} for the Penrose Inequality, in the sense that, arbitrarily close to $(M,g)$ in $C^{2,\alpha}$, we can simultaneously decrease the mass of the end and increase the area of the outermost minimizing surface:

\begin{restatable}{theorem}{nopenroseqf}\label{thm: no_penrose_qf}
     Let $(M,g)$ be the one ended hyperbolic manifold obtained by cutting a quasi-Fuchsian manifold (that is not Fuchsian) along a strictly stable outermost minimal surface, and considering an end that contains no compact minimal surfaces. For any $\e > 0$, there is a $(\bar{M},\bar{g})$ with boundary, ALH at rate $q = 3$ to $(M,g)$ with the following properties:
     \begin{enumerate}
     \item The scalar curvature of $\bar{g}$ satisfies $R(\bar{g})= -6$,
         \item The boundary $\del \bar{M}$ is outer-minimizing and minimal,
         \item The area of the boundaries satisfy $\text{Area}(\del \bar{M},\bar{g}) > \text{Area}(\del M,g)$,
         \item The mass of $(\bar{M},\bar{g})$ is negative, and
         \item The norm $|\bar{g}-g|_{C^{2,\alpha}}< \e$.
     \end{enumerate}
\end{restatable}

This document is organized as follows. In Section \ref{sect:prelim} we record definitions and recall key facts about static manifolds. In Section \ref{sect: static potentials are proper}, we prove Theorem \ref{thm: level sur are closed}. In Section \ref{sect: doubles}, we prove the uniqueness results: Theorems \ref{thm: bonus} and \ref{thm:main}. In Section \ref{sect: quasi}, we use the results from Section \ref{sect: doubles} to prove that there is a family of one ended hyperbolic manifolds which are not optimal for the Penrose inequality. \\

\noindent {\bf Computational Resource Disclosure.} LLMs were used to search the literature and check proofs. An LLM suggested a known technique to prove Theorem \ref{thm: bonus}.\\

\noindent {\bf Acknowledgments. } I am grateful to the continued support of my PhD advisor, Andr\'e Neves. I am grateful to Lan-Hsuan Huang for explaining the state of the literature to me.

\section{Preliminaries}\label{sect:prelim}
In this section, we discuss mean curvature convention, give basic definitions, review fundamental properties of static manifolds, and recall known results relevant to this paper. \\

\noindent {\bf Conventions.} We define the second fundamental form of a surface as \[A(X,Y) = \inner{\grad_X \nu}{Y}\] where $X,Y$ are tangent to $S$ and $\nu$ is the upward pointing unit normal. Then \[H = g^{ij}A_{ij}\] and \[\vec{H} = -H\nu.\] 

\subsection{Asymptotically Locally Hyperbolic (ALH) Manifolds} 
We define and discuss the basic objects of study in this paper, asymptotically locally hyperbolic (ALH) manifolds.

\begin{definition}\label{def: reference end}
    Let $k\in\{0,\pm 1\}$ and let $N_k$ is an $n-1$ dimensional closed manifold quipped with an Einstein metric $b$ such that $\ric(b) = k(n-2)b$. Let $\model = (1,\infty) \times N_k$ and \[\mg = \frac{dr^2}{r^2+k} + r^2b.\]  
    We refer to the pair $(\model,\hat{g})$ as a \emph{model space} or \emph{reference end}. 
\end{definition}
Let $\{\del_{x_i}\}_{1\leq i\leq n-1}$ be an orthonormal basis for $(N_k,b)$. We sometimes use the following orthonormal basis for $(\model,\mg)$: \[e_r = \sqrt{r^2+k}\;\;\; \del_r, e_{x_i} = r^{-1}\del_{x_i},\] for $1\leq i\leq n-1$. 
The mean curvature of the $r$- level surfaces of $(\model,\mg)$ is \[H(r) = (n-1)\frac{\sqrt{r^2+k}}{r}.\] A useful change of coordinates is given by 
\begin{itemize}
    \item If $k = 1$, set $\sinh(t) = r$. The metric becomes \[\mg = dt^2 + \sinh^2(t)b.\] 
    \item If $k = 0$, set $r = e^t$. The metric becomes \[\mg = dt^2 + e^{2t}b.\]
    \item If $k=-1$, set $r = \cosh(t)$. The metric becomes \[\mg = dt^2 + \cosh^2(t)b.\] In this case, the boundary $r = 1$ or equivalently $t =0$ is totally geodesic. We can double $\model$ over its boundary to obtain a complete and in fact analytic metric  $dt^2 + \cosh^2(t)b$ on $\R\times N$. 
\end{itemize}

We next introduce the function spaces we will use throughout the paper. Let $B\subset\model$ be compact and let $\hat{\grad}$ be the connection on $\model$:

\begin{definition} (\cite{huang22}, Definition 2.5)
    Let $k \in \N_{\geq 0},\alpha\in (0,1)$, and $q\in\R$. The \emph{weighted H\"{o}lder spaces} $\ckaq (\model \setminus B)$ are the collection of $C^{k,\alpha}_{\text{loc}}$ functions that satisfy \[\|f\|_{\ckaq(\model\setminus B)} := \sum_{|I| = 0,1,...k}\;\sup_{x\in\model\setminus B} |x|^q|\moddec{\grad}^{I}f(x)|_{\mg} + \sup_{x\in\model\setminus B} |x|^{q}[\moddec{\grad}^kf]_{\alpha;B_1(x)} < \infty\]
    and 
    \[[{\grad}^kf]_{\alpha;B_1(x)} = \sup_{1\leq i_1,...,i_k\leq n, \; y\neq z\in B_1(x)}\frac{|e_{i_1}\cdots e_{i_k}f(y) - e_{i_1}\cdots e_{i_k}f(z)|}{d(y,z)^{\alpha}},\] where $B_1(x)$ is the unit ball centered at $x$ intersecting $\model\setminus B$, $\{e_i\}_{1\leq i\leq n}$ is the orthonormal frame described in Definition \ref{def: reference end} for $\mg$ and $|x| = r$
    \end{definition}
    
    This definition extends to tensors by saying a tensor $h\in  \ckaq(\model\setminus B)$ if and only if the components of the tensor with respect to an orthonormal frame lie in $\ckaq(\model\setminus B)$.

    To formally define weighted H\"{o}lder spaces for general ALH manifolds, let $M$ be a manifold covered by an atlas that consists of a non-compact chart: $\Phi\colon \model\setminus B\to M\setminus K $, for compact sets $B$ and $K$, and a finite number of compact charts. The \emph{weighted H\"older norm} $\|f\|_{C^{k,\alpha}_{-q}(M)}$ (for functions and tensors) is the sum of $\|\Phi_*f\|_{C^{k,\alpha}_{-q}(\model\setminus B)}$ and the usual $C^{k,\alpha}$ norms on compact charts. 
    
    By $C^{k,\alpha}_{-q}(M)$ we denote functions in $C^{k,\alpha} (M)$ with finite weighted H\"older norm. Additionally, we say that a function is $O(r^{-q})$ if it belongs to $C^{2,\alpha} (M)$ and we say a tensor is $O(r^{-q})$ if all its components belong to $C^{2,\alpha}(M)$ when considered with respect to the orthonormal basis. 

With this function space definition in mind, we may formally define the main objects of study in this paper. 

\begin{definition} (\cite{huang22} Definition 2.6)
\label{def:alh}
    Let $M^n$ be a  $n$-dimensional connected manifold with smooth Riemannian metric $g$. We say that $(M,g)$ is \emph{asymptotically locally hyperbolic (ALH)} with respect to a reference end $(\model,\mg)$ at rate $q>0$ if the following hold:
    \begin{enumerate}
        \item There exists a diffeomorphism $\Phi \colon M\setminus K \to \model\setminus B$ where $K\subset M$ and $B\subset \model$ are compact sets. The induced coordinate chart is called the \emph{chart at infinity}.
        \item With respect to the chart at infinity, $g - \mg\in O_2(r^{-q})$.
        \item The scalar curvature of $g$, $R(g)$, satisfies $r(R(g) + n(n-1)) \in L^1(M\setminus K)$
    \end{enumerate}
\end{definition}
We will often refer to a connected component of $M\setminus K$ as an \emph{exterior region} or \emph{asymptotic region} of the manifold, which we denote by $M_{ext}$. Since $\Phi$ gives a diffeomorphism from $M_{ext}$ to an end of $\model$, $M_{ext}$ has coordinates $(r,x)$ induced by $\Phi$. If all ends of a manifold $(M,g)$ are ALH at rate $q$, we say $(M,g)$ is ALH at rate $q$.  

To more easily discuss the geometry of ALH manifolds $M$, we make the following definition:

\begin{definition}\label{def: conf inf}
    Let an end of $M$ be ALH to a particular reference end $\model = (1,\infty)\times N, \hat{g} = \frac{dr^2}{r^2+k} + r^2b$. We say that $(N,b)$ is the \emph{conformal infinity} of $M$. We say that the conformal infinity of $M$ has positive, flat, or negative type, corresponding to the curvature of $(N,b)$. 
\end{definition}

\subsection{Static manifolds} 
In this section, we discuss known results about static manifolds, which are manifolds that admit solutions to a certain overdetermined system of PDE. 

\begin{definition}\label{def: static}
    A Riemannian manifold $(M^n,g)$ (possibly with boundary) is \emph{static} if it admits a nontrivial solution to \[L_g^*(f) =-(\lap f )g + \hess_g f -f\ric(g) = 0. \] The function $f$ is called a \emph{static potential} on $(M,g)$. If $\del M\neq \emptyset$, we require $ f = 0$ on $\del M$ and in that case we call $(M,g)$ a \emph{static black hole}. 
\end{definition}
We sometimes refer to the triple $(M,g,f)$ as a \emph{static triple}. The following class of examples, as discussed in the introduction, is the direct generalization of Schwarzschild metrics to the asymptotically (locally) hyperbolic setting.

\begin{example}\label{ex: kottler} (Kottler metrics) Let $(N,b)$ be a closed manifold with $\ric(b) = k(n-2)b$ for $k\in\{0,\pm 1\}$. Let $V(r) =\sqrt{r^2+k-\frac{2m}{r^{n-2}}}$ and let $m\in \mathbb{R}$ be large enough so that $V(r)^2$ has a simple, positive zero, and let $r_m$ be the largest zero of $V(r)$. Define \[g = \frac{dr^2}{V(r)^2} + r^2b\] on $M = (r_m,\infty)\times N^{n-1}$. The boundary $r = r_m$ is a totally geodesic hypersurface and $(M,g,V(r))$ is a static ALH black hole. 
\end{example}

The next group of lemmas record some well known but important facts that we will use about static manifolds.

\begin{lemma}\label{lem: static means constant scalar curv}
    Let $(M,g)$ be a static manifold. The scalar curvature of $M$ is constant. 
\end{lemma}
\begin{proof}
    See, for example (\cite{corvino2000}, Proposition 2.3).
\end{proof}

\begin{lemma}\label{lem:fisconvex}
   For large enough $r$, the potential $f$ is strictly convex when $f > 0$ and strictly concave when $f < 0$.  
\end{lemma}
\begin{proof}
    If $g$ is ALH at rate $q$, $\ric(g) = -(n-1)g + O(r^{-q}).$ So, 
\begin{align*}
    \hess f &= f(\ric(g) + ng) \\
    &= f(g + O(r^{-q})),
\end{align*}
meaning that $\hess(f)$ is definite when $r$ is large enough.
\end{proof}
\begin{lemma}\label{lem: general prop of statics}
    Let $(M^n,g)$ be an ALH static manifold of order $q$ with nontrivial static potential $f$. Let $\s_s = \{p\in M | f(p)=s\}$. Then,
    \begin{enumerate}
        \item  When $|\grad f|\neq 0,$ \[H(\s_s) = \frac{s}{|\grad f|}(n-1 + O(r^{-q})),\] where $\grad f$ is being evaluated on $\s_s$. 
        \item If $\s_0$ is nonempty, then it is an embedded, separating totally geodesic hypersurface. 
    \end{enumerate}
\end{lemma}
\begin{proof}
    Let $\del_i,\del_j$ be basis vectors for $\s_s$ and let $A(\s_s)$ be the second fundamental form of $\s_s$. We first calculate: 
    \[
        A(\s_s)_{ij} = \inner{\grad_i\frac{\grad f}{|\grad f|}}{\del_j}
        = \frac{\hess f(\del_i,\del_j)}{|\grad f|}
        = \frac{s(\ric(g) + ng)_{ij}}{|\grad f|}, \tag{$\star$}
\]
when $|\grad f|\neq 0$. 

\noindent (1): When $|\grad f|\neq 0,$   plugging in $\ric(g) + ng = g + O(r^{-q}),$  and taking the trace of the rightmost equality in ($\star$) above and combining with Lemma \ref{lem: static means constant scalar curv} implies \[H(\s_s) = \frac{s}{|\grad f|}(n-1 + O(r^{-q}))\] as desired. 

\noindent (2): We first show that 0 isn't a critical value of $f$ when $f$ is not identically 0. Assume there is some $p\in\S_0$ such that $\grad f(p) = 0.$ 
If $\gamma$ is any unit speed geodesic crossing through $p$ at $t =0$, by the static equations there is a function $a(t)$ so that $f^{\p\p}(\gamma(t)) = a(t)f(\gamma(t))$ for small enough $|t| < \e $ and $f(\gamma(0)) = f^{\p}(\gamma(0)) = 0$. By uniqueness of solutions to ODEs, $f\equiv 0$ along $\gamma$. Since $M$ is connected, this implies that $f = 0$ on all of $M$, which is a contradiction. 

Therefore, $\s_0$ is embedded, and it is separating since it is the level set of a function. 

\end{proof}

\begin{lemma}\label{lem: surface gravity}
    Let $f$ be a static potential, let $\s_0$ be the zero set of $f$. Then, $|\grad f|$ is constant on each connected component of $\s_0$. 
\end{lemma}
\begin{proof}
    Let $\del_i$ be tangent to $\s$. Then
    \begin{align*}
        \del_i|\grad f|^2 &= 2\hess(f)(\del_i,\grad f)
    \end{align*}
    Since $\hess (f) = 0$ on $\s$, $\del_i|\grad f|^2 = 0$.
\end{proof}

Lastly, we record a useful theorem. 
\begin{theorem}(\cite{Ch05}, Section 4)\label{thm: doubled is analytic}
    Let $(M,g)$ be a static black hole with static potential $f$. The doubling of $(M,g)$ across its black hole boundary is globally analytic and the odd reflection of the potential is an analytic function.
\end{theorem}

\section{Asymptotic Behavior of Static Potentials}\label{sect: static potentials are proper}
In this section, we prove Theorem \ref{thm: bonus}.

\begin{lemma}\label{lem: uniform bounding of H for static}
    Let $f$ be a static potential such that $f > 0$ in $M_{ext}$. Near infinity, there exists a $c>0$ such that $|\grad f| < cf$. 
\end{lemma}

\begin{proof}
    Let $\{x_i\}\subset M_{ext}$ be a sequence of points such that $|x_i|\to\infty$. Let $B_{x_i}(r)$ be the geodesic ball of radius $r$ centered at $x_i$. Fix some $\e > 0$ and let $a_i = \dist(x_i,\del M_{ext}) - \e$ and assume that $x_i$ is chosen so that $a_i > 1$ for all $i$. This way, $B_{x_i}(a_i)\subset M_{ext}$. 
    
    Then, (\cite{ChYa75}, Theorem 6) implies that there exist $\alpha_n>0$  depending only on dimension and a $K$ depending only on $g$ such that \[|\grad f(x_i)|\leq \alpha_nf(x_i)(|K| + \frac{1}{a_i}).\] Setting $c = \alpha_n(|K|+1)$ gives the desired conclusion. 
\end{proof}

 Our next lemma guarantees the existence of a geodesic with certain properties that will be used to create a contradiction to the maximum principle in the main proof.

\begin{lemma}\label{lem:geodesics contained in cones} Let $(M^n,g)$ be an ALH static manifold with a static potential $f$ such that $f > 0$ near infinity. If $x_1$ and $x_2$ are points on the conformal boundary of $(M,g)$, there is a curve $\sigma\colon(-\infty,\infty)$ connecting them such that the potential $f$ is convex along $\sigma$ if $x_1$ and $x_2$ are sufficiently close together. 
 
\end{lemma}
\begin{proof}  First, let $k\in\{0,\pm 1\}$. Let $N$ be an $n-1$ dimensional manifold that is equipped with a metric $b$ where $\ric(b) = k(n-2)b$. Let $M = (0,\infty)\times N$. Set $t = \log(r)$ and consider the metric  \[\tilde{g} = dt^2 + e^{2t}b.\] We construct a geodesic connecting any two points $x_1,x_2$ on the conformal boundary of $(M,g)$.  Note that in any case, for any $\tilde{g}$, there is a model space $g$ such that $|\tilde{g}-g| = O(r^{-2})$. Therefore, for an arbitrary ALH static metric (at rate $q$), there is a $\tilde{g}$ such that \[|\tilde{g}-g| = O(r^{-a}),\] where $a = \min\{2,q\}$. 

Let $\gamma\colon (-\infty,\infty)\to \tilde{M}$ be unit-speed geodesic. We can write $\gamma(s) = (t(s),x(s))$. The parametrization by arc length implies 
\begin{equation}\label{eq: pbal}
    (t^{\p}(s))^2 + e^{2t(s)}|x^{\p}(s)|_b^2 = 1.
\end{equation} Combining with the geodesic equation for $t^{\p\p}$ gives 
\[t^{\p\p}(s) = 1-(t^{\p}(s))^2,\] 
which, choosing a non-radial solution gives 
\[t(s) = \ln(\cosh(s-C_0))+C_1,\] 
for two constants $C_0$ and $C_1$. We assume that $C_0 = 0$ so that the minimum $t$ value occurs at $t = C_1$. In turn, plugging the equation for $t(s)$ back into (\ref{eq: pbal}) gives $|x^{\p}(s)|_b = e^{-C_1}\text{sech}^2(s)$. 

Let $(N,b)$ be the conformal infinity of $\tilde{M}$. Let $x_1,x_2\in N$, $x_1\neq x_2$ and let $\ell = \dist_b(x_1,x_2)>0$. Let \[\tilde{\gamma}\colon \Big[-\frac{\ell}{2},\frac{\ell}{2}\Big]\to N\] be a length minimizing geodesic in $(N,b)$, parametrized by arc length such that $\tilde{\gamma}(-\frac{\ell}{2}) = x_1$ and $\tilde{\gamma}(\frac{\ell}{2}) = x_2$. Then, if $x(s) = \tilde{\gamma}(\varphi(s))$, we see $\varphi^{\p}(s) = e^{-C_1}\text{sech}^2(s)$ so \[x(s) = \tilde{\gamma}(e^{-C_1}\tanh(s)).\] Picking $e^{-C_1} = \frac{\ell}{2}$ shows that the desired geodesic connecting the two points $x_1,x_2$ on the conformal boundary is given by \[\gamma(s) = \left(\ln(\cosh(s)) + \ln\frac{2}{\ell},\tilde{\gamma}(\frac{\ell}{2}\tanh(s))\right).\] In particular, as $\ell\to 0,$ $\ln\frac{2}{\ell}\to\infty$, showing that the minimum height of the geodesic goes to infinity as $x_1\to x_2$.

Let $a = \min\{q,2\}$. For an arbitrary ALH manifold $(M,g)$, consider the space $\tilde{g}$ so that $|\tilde{g}-g|_{g} = O(r^{-a})$. Let $x_1$ and $x_2$ be points on the conformal boundary and let $\sigma$ be the geodesic connecting $x_1$ and $x_2$ for $\tilde{g}$. We show that if $(M,g)$ is static, then a static potential $f$ for $g$ is convex along $\sigma$ as long as $x_1$ and $x_2$ are sufficiently close. 

Next, for an arbitrary ALH manifold $(M,g)$, consider the space $\tilde{g}$ so that $|\tilde{g}-g|_{g} = O(r^{-a})$.  Let $\sigma$ be a geodesic for $\tilde{g}$. Let $\varphi$ be a function. We have \[\frac{d^2}{dt^2}\varphi(\sigma(t)) = \hess_{g}\varphi(\sigma^{\p},\sigma^{\p})  +d\varphi((\grad-\tilde{\grad})_{\sigma^{\p}}\sigma^{\p})\] where $|\tilde{\grad}-\grad|_g = O(r^{-a})$. 

Let $\varphi$ be a static potential $f$ for the metric $g$. Since $|\grad f|_g<cf$ by Lemma \ref{lem: uniform bounding of H for static}, \[|df((\grad-\tilde{\grad})_{\sigma^{\p}}\sigma^{\p})|\leq |\grad f||(\tilde{\grad}-\grad)_{\sigma^{\p}}\sigma^{\p}|< Cfr^{-a}.\]

Note that $\inner{\sigma^{\p}}{\sigma^{\p}}_g = 1+ O(r^{-a})$. Thus, \begin{align*}
    \hess_gf(\sigma^{\p},\sigma^{\p}) &= f(\ric(g) + ng)(\sigma^{\p},\sigma^{\p})\\
    &= f(1 + O(r^{-a})).
\end{align*}

Putting it all together, 

\begin{align*}
    \frac{d^2}{dt^2}f(\sigma(t)) &= \hess_{\tilde{g}}f(\sigma^{\p},\sigma^{\p}) \\
    &= \hess_{\tilde{g}}f(\sigma^{\p},\sigma^{\p}) + df((\tilde{\grad}-\grad)_{\sigma^{\p}}\sigma^{\p})\\
    &> f(1+O(r^{-a})) - Cfr^{-a}\\
    &=f(1+O(r^{-a})),
\end{align*}
which is positive once $r$ becomes sufficiently large. 
\end{proof}

   \begin{lemma}\label{lem: curves at finite distance}
       Let $\sigma_1,\sigma_2$ be two curves such that $\lim_{t\to\infty}|\sigma_1(t)| = \lim_{t\to\infty}|\sigma_2(t)| = \infty$ and $\sigma_1$ and $\sigma_2$ remain at finite distance from each other near infinity. Then, if $f$ is a static potential on $M$ such that $f\geq 0$ and as $t\to\infty$, 
       \begin{enumerate}
           \item $f(\sigma_1(t))\to\infty$, or
           \item $f(\sigma_1(t)) \to 0$,
       \end{enumerate}
       then $f(\sigma_2(t))$ has the same limit as $t\to\infty$. 
   \end{lemma}

\begin{proof}
    Let $\sigma_1$ and $\sigma_2$ be at distance less than $d$ from each other near infinity. By Lemma \ref{lem: uniform bounding of H for static}, there is a constant $c>0$ such that $f>c|\grad f|$. Let $t_0\geq0$ and let $T$ be the time such that $\sigma_1(T)$ is a closest point on $\sigma_1$ to $\sigma_2(t_0)$. Let $\bar{\sigma}\colon[0,\ell]\to M$ be the geodesic parameterized by arc length that connects $\sigma_1(T)$ and $\sigma_2(t_0)$. The fundamental theorem of calculus applied to $\log(f)$ then implies that $f^{\p}(\bar\sigma(t))\leq cf(\bar\sigma(t))$, meaning that $f(\sigma(t))e^{-cd}\leq f(\sigma_2(t_0))\leq f(\sigma_1(T))e^{cd}$. If $f(\sigma_1(t))\to 0$, taking $t\to\infty$ gives the desired result.   
\end{proof}
Finally, we collect some facts about the asymptotics of the static potential. We call a geodesic ray $\gamma\colon [0,\infty)\to M$ a \emph{radial geodesic} if $\gamma$ is a geodesic emanating from $\del M_{ext}$ with velocity $\del_r/|\del_r|$.

\begin{lemma}(\cite{HuJaMa}, Theorem 3.5) \label{lem: asymptotic properties of f along geodesics}
    Let $f$ be a static potential on an ALH manifold $M$.
Let $\gamma(t)$ be a radial geodesic ray. Assume $f\geq 0$. There are two possible end behaviors: 
        \begin{enumerate}
            \item As $t\to\infty$, $f(\gamma(t))\to\infty$ or
            \item As $t\to\infty$, $f(\gamma(t))\to 0$.
        \end{enumerate}
   \end{lemma}

\begin{lemma}\label{lem: uniformity}
       Let $\gamma(t)$ be a geodesic ray. If $f(\gamma(t))\to\infty$ as $t\to\infty$, then there is an open ball $B\subset N$ containing $\gamma(0)$ such that there is a constant $C> 0$ such that if $x\in B$, \[C^{-1}r\leq |f(x,r)|\leq Cr.\]  Further, if for all geodesic rays $\eta$ emanating from a compact $K\subset\del M_{ext}$ have $f(\eta(t))\to\infty$ as $t\to\infty$, then there exists a constant $\bar{C}>0$ such that \[\bar{C}^{-1}r\leq |f(x,r)|\leq \bar{C}r\] for all $x\in K$.
\end{lemma}
   \begin{proof}
    The first conclusion is (\cite{HuJaMa}, Theorem 3.5). We move next to uniformity. For each $c > 0$, we consider $D(x,c)$ to be largest open set containing $x\in \del M_{ext}$ such that \[c^{-1}r< |f(x,r)|< cr.\] The set $D(x,c)$ may be empty. A cover of $K$ is defined by 
    \[ K \subset \bigcup_{x\in M}\bigcup_{c>0} D(x,c).\]
 Since $K$ is compact, this admits a finite open sub-cover: \[K \subset D(x_1,c_1)\cup\cdots\cup D(x_n,c_n)\] for some $n\in \N_{>0}$. Then letting \[\bar{C} = \max_{0\leq i \leq n} c_i\] shows that \[\bar{C}^{-1}r< |f(x,r)|< \bar{C}r\] for all points $(x,r), x\in K$. 
   \end{proof}

   \begin{lemma}\label{lem: finite dist}
Let $\sigma\colon(-\infty,\infty)\to M$ be the curve that connects the points $x_1$ and $x_2$ at infinity along which $f$ is convex constructed in Lemma \ref{lem:geodesics contained in cones}. Then $\sigma$ remains at finite distance from the geodesic ray $\gamma_{x_1}$ that intersects the conformal boundary at $x_1$ as $t\to-\infty$ and $\sigma$ remains at finite distance from the geodesic ray $\gamma_{x_2}$ that intersects the conformal boundary at $x_2$ as $t\to\infty$
   \end{lemma}
\begin{proof}
We turn our attention to showing that as $t\to-\infty$, $\sigma$ remains at finite distance from the geodesic ray $\gamma_{x_1}$ that intersects the conformal boundary at $x_1$. The case for $\sigma$ staying at finite distance from $\gamma_{x_2}$ near infinity works analogously. 

Let the coordinate ray $r(s,x_0)$ be defined as\[r(s,x_0) =(s,x_0).\] Next, define $\eta_s\colon[0,1]\to M$ by \[\eta_s(t) = (\ln(\cosh(s)) +\ln (2/\ell),\tilde{\gamma}(-\frac{t\cdot \ell}{2} + (1-t) \frac{\ell}{2}\tanh(s)) ).\] Note that \[\eta_s(0) = (\ln(\cosh(s)) +\ln (2/\ell), \tilde{\gamma}(\frac{\ell}{2}\tanh(s))),\] which is a point on $\sigma$. And, \[\eta_s(1) = (\ln(\cosh(s)) +\ln (2/\ell), x_1 ),\] which is on the coordinate ray $r(t,x_1)$.

By direct computation, the length of $\eta_s$ is less than $Ce^s$ where $C$ is a constant that depends on $g$, but not $s$. Therefore, as $s\to-\infty$, the length of this curve goes to 0.

Next, we show that the coordinate ray $r(t,x_1)$ and $\gamma_{x_1} = (t(s),x(s))$ remain at finite distance from each other. Since $\tilde{g}$ is ALH to $g$, there is some constant $C$ such that:
\[t^{\p}(s)^2 + e^{2t(s)}|x^{\prime}(s)|^2 = |\gamma_{x_1}^{\p}(s)|_{\tilde{g}} \leq C. \] Therefore, $|x^{\prime}(s)|_b\leq \sqrt{C}e^{-t(s)}$. As in (\cite{HuJaMa} line 3.13), $t(s) = s+O(1)$. And thus, \[\int_{s}^{\infty}|x^{\p}(u)|du\leq Ce^{-s}, \] as $s\to\infty$, $x(s)$ converges to some $x_{\infty}$. Then, \begin{align*}
    d_g((t(s),x(s)),(t(s),x_{\infty}))&\leq Ce^{t(s)}d_b(x(s),x_{\infty})\leq C^{\p}.
\end{align*}
Therefore, if $x_1 = x_{\infty}$, $r(t,x_1)$ and $\gamma_{x_1}$ remain at finite distance from each other near infinity. 
\end{proof}
   
We combine our efforts to arrive at the main result of this subsection. 
\proper*

\begin{proof}
Let $\gamma_p$ be a geodesic ray from $p\in\del M_{ext}$. Let \[ F = \{x\in \del M_{ext} :
\lim_{t\to\infty}f(\gamma_x(t))=0\}.
\] 
We show that distinct $p\in F$ are separated by some minimum distance $\delta$.
Suppose for the sake of a contradiction that there exist points $p_i\in F$ such that for all $i\in\N$, we can find $x_i\in F$ so that $d_b(x_i,p_i)<\frac{1}{i}$, where $d$ is the distance from $g$. 
Let $\gamma_{x_i}$ and $\gamma_{p_i}$ be the radial geodesics intersecting the conformal boundary at $x_i$ and $p_i$ respectively. 
By Lemma \ref{lem:geodesics contained in cones}, there is a curve $\sigma_i\colon (-\infty,\infty)\to M$ that intersects the conformal boundary at $x_i$ and $p_i$ such that the minimum $r$ coordinate of $\sigma_i$ goes to $\infty$ as $\dist(x_i,p_i)\to 0$. By Lemma \ref{lem: finite dist} $\sigma_i$ remains at finite distance from $\gamma_{x_i}$ near $t = -\infty$ and near $\gamma_{p_i}$ near $t = \infty$ and so Lemma \ref{lem: curves at finite distance} and Lemma \ref{lem: asymptotic properties of f along geodesics} combined imply that $\lim_{t\to-\infty}f(\sigma_i(t)) = \lim_{t\to\infty}f(\sigma_i(t) )= 0$. Since $f>0$, $f( \sigma_i)$ must have an interior maximum. Taking $i$ large enough guarantees that $f(\sigma_i)$ cannot have an interior maximum, which creates a contradiction.

Since all points in $F$ are separated from all other points by $\delta$-neighborhoods, the compactness of $N$ implies that $F$ is a finite set. 

The second conclusion of Theorem \ref{thm: level sur are closed} follows from Lemma \ref{lem: uniformity}
\end{proof}

\section{Uniqueness of Analytic Doubles}\label{sect: doubles}

Let $M$ be a manifold with boundary. The \emph{double} of $M$, $D(M)$, is defined as \[D(M) =M\times \{0,1\}\setminus \sim,\] where $(x,0)\sim (x,1)$ for all $x\in \del M$. Throughout this section, we use the coordinates $(x,t)$ whose properties will be specified below. To clarify notation: by $g(x,t)$ we mean the full $n\times n$ metric tensor of $g$ evaluated at the point $(x,t)$. By $\del_i g(x,t)$ we mean the $i$th partial derivative of $g(x,t)$ and by $\del_i g(x,0)$ we mean the $i$th partial derivative of $g(x,t)$ evaluated at the surface $t = 0$. 

\bonus*

\begin{proof}

We consider the double of $M$ across its boundary, $D(M)$. Note that by Theorem \ref{thm: doubled is analytic}, $D(M)$, along with the doubled metric, is analytic. In a neighborhood of $\s_0\subset D(M)$, we can write $g$ in terms of an equidistant foliation from $\s_0$:\[g = dt^2 + h_t(x),\] where $h_t$ is a metric on the surface of signed distance $t$ from $\s$. Note that this is an analytic choice of coordinates because the exponential map is analytic.  This coordinate system has the property that $g(x,0)_{tt} = 1, g(x,0)_{ti} = 0$, and $\del_t^m g(x,0)_{tt} = \del_t^mg(x,0)_{ti} = 0$ when $m\geq 1.$

We now use the static equations to show that all the derivatives of a pair $(f,g)$ are determined on $\s_0$, which will give uniqueness essentially because $f$ and $g$ are analytic. 

 By Lemma \ref{lem: surface gravity}, we can normalize so that $\del_tf(x,0) = 1$. Consider the following Cauchy data:
    \begin{align*}
        (g(x,0),f(x,0)) &= (dt^2+g_{\s_0},0) \\
        (\del_t g(x,0),\del_tf(x,0)) &= (0,1) \\
    \end{align*}

Now, we will use induction to show that all derivatives of order $m$ of a preexisting solution $(f,g)$ are determined by the Cauchy data on $\s_0$. 

If $m = 1$, all derivatives of order $m$ are part of the Cauchy data. We now assume that the derivatives of order up to $m$ are determined by the Cauchy data and show that the derivatives of order $m+1$ are determined by the Cauchy data. 

For convenience, by $F_m(x)$, we denote a function of the derivatives $\del_{x}^\alpha\del_t^l(f,g)$, where $x$ represents tangential derivatives and $\alpha$ is a multi-index of arbitrary magnitude, and $0\leq l\leq m$. To keep notation clean, we allow $F_m(x)$ to change line to line. 

Therefore, if all derivatives up to order $m$ are determined by the Cauchy data, then all functions expressible as $F_m(x)$ are also determined by the Cauchy data. Therefore, we are left to show that the derivatives of order $m+1$ are expressible as $F_m(x)$. Therefore, the only derivatives of order $m+1$ that must be checked explicitly are $\del_t^{m+1}f$ and $\del_t^{m+1}g$.  

First, we show $\del_t^{m+1}f = F_m(x)$. Indeed,
\begin{align*}
    0 = \del_t^{m-1}(\lap_g f -nf)&= \del_t^{m-1}(g^{ij}(\del_i\del_jf-\Gamma_{ij}^k\del_kf)-nf)\\
    &= \del_t^{m+1}f + F_m(x).
\end{align*}

Next, we show that $\del_t^{m+1}g = F_m(x)$. We first recall from our choice of coordinates, $\del^{m+1}_{t}g(x,0)_{tt} = \del_t^{m+1}g(x,0)_{ti} = 0$. Therefore, we are left to handle the cases of $\del_t^{m+1}g(x,0)_{ij}$ where $i,j\neq t$. 

We will differentiate the equation $\hess f = f(\ric(g) + ng)$. On one hand,
\begin{align*}
    \del_t^m\hess f_{ij} &= \del_t^m (\del_i\del_j f -\Gamma_{ij}^k\del_kf)\\
    &= \del_t^m \del_i\del_j f - \del_t^m (\Gamma_{ij}^k\del_kf ) 
\end{align*}
where \[\del_t^m (\Gamma_{ij}^k\del_kf ) = \del_t^m(\Gamma_{ij}^k)\del_kf +F_m(x).\] 
Now, on $\s_0$, $\del_i f = 0$ if $i\neq t$, and $\del_t f = 1$.  Since if $i,j\neq t$, $\del_t^m\del_i\del_j f = F_m(x)$ on $\s_0$.
And, because the cross metric terms vanish to infinite order on $\s_0$, \[\del_t^m(\Gamma_{ij}^t) =-\frac{1}{2}\del_t^mg_{ij,t} .\]Therefore, 
\[\del_t^m\hess f = \frac{1}{2}\del_t^{m+1}g_{ij} + F_m(x).\] on $\s_0$. 

Now we look at $\del_t^m(f(\ric(g) +ng))$. Because $f$ vanishes on $\s_0$, the highest order derivative term for $\ric(g) + ng$ is the term $m\del_t f \del_t^{m-1}(\ric(g) +ng)$. Since $n\del_t^{m-1}g = F_m(x)$, we compute $\del_t^{m-1}\ric(g)_{ij}$.

Let $A$ be the second fundamental form of $\s_0$ and let $H$ be its mean curvature. Let $\ric_M$ be the Ricci tensor for $M$ and let $\ric_{\s_0}$ be the Ricci tensor for $\s_0$. Let $R_M$ be the Riemannian curvature tensor for $M$. The Gauss-Codazzi equations give:
\begin{align*}
    \ric_M(\del_i,\del_j) &= \ric_{\s_0}(\del_i,\del_j) + \inner{R_M(\del_i,\del_t)\del_t}{\del_j} + A_{i}^kA_{kj} - HA_{ij}.
\end{align*}
Therefore, 
\begin{align*}
    \del_t^{m-1}\ric_M(\del_i,\del_j) &= \del_t^{m-1}\ric_{\s_0}(\del_i,\del_j) + \del_t^{m-1}\inner{R_M(\del_i,\del_t)\del_t}{\del_j} +\del_t^{m-1}( A_{ik}A_{kj} - HA_{ij}).
\end{align*}

We examine the behavior of each of the three terms independently. First, since $\ric_{\s_0}$ involves only tangential derivatives of $g$, $\del_t^{m-1}\ric_M(\del_i,\del_j) = F_m(x)$. 

Note that in coordinates, $(R_M)_{ittj}$ is:
\[R_{ittj} = -\frac{1}{2}  \partial_{t}\partial_{t} g_{ij}  + F_1(x)\] 
The term $\del_t^2 g_{ij}$ is a second order derivative in time, while all other terms contain at most first order derivatives in time. Therefore, 

So, \[\del_t^{m-1}(R_M)_{ittj} = -\frac{1}{2} \del_t^{m+1}g_{ij} +F_m(x).\]

Since $\s_0$ is totally geodesic, and since $H$ and $A$ depend on the first derivatives of $g$, \[\del_t^{m-1}( A_{i}^kA_{kj} - HA_{ij}) = F_m(x).\]

Putting both sides of the equation back together, we see that \[\frac{1}{2}\del_t^{m+1}g_{ij} = -\frac{m}{2}\del_t^{m+1}g_{ij} + F_m(x),\] which means that $\del_t^{m+1}g$ is entirely determined by the Cauchy data. Therefore, since $g$ is analytic, $g$ is uniquely determined in a neighborhood of $\s\subset D(M)$. Since $\s$ is separating and $\pi_1(\s)=\pi_1(M)$, this means that $g$ is uniquely determined on a collar neighborhood of $\s_0$. Then, [Theorem 5.2, \cite{KoWa24}] implies that $g$ is actually uniquely determined on all of $D(M)$, which in turn implies that $g$ is the unique metric on $M$ such that $g_N  = g_{\bh}$. 
\end{proof}

\main*

\begin{proof}
    We suppress $\genus\geq 2$ in notation and refer to $\s_{\genus\geq 2}$ simply as $\s$. Let $\Gamma$ be the set of surfaces isotopic to $\{0\}\times\s$. Let \[A = \inf_{S\in\Gamma} \text{Area}(S).\] Let $S_i\in \Gamma$ be such that $\text{Area}(S_i)< A + 1/i$. From the ALH assumption, near either end there is a foliation by mean convex surfaces coming from the coordinate foliation in the model, meaning that $S_i$ can be taken to all be contained in a bounded region. As a consequence of (\cite{MeSiYa82}, Theorem 1), there is an embedded least area minimal surface $S_{\infty}$ isotopic to $\s$. 
    
    We now claim that $ f= 0$ along $S_{\infty}$. Following the arguments in the proof of (\cite{huang22}, Lemma 5.4), if $f$ is not identically $0$ along $S_{\infty}$, then $f$ has a sign along $S_{\infty}$ and a neighborhood of $S_{\infty}$ is foliated by totally geodesic hypersurfaces. Replacing $f$ with $-f$ if necessary, we may assume $f > 0$ on $S_{\infty}$. However, this is not possible. A foliation by totally geodesic hypersurfaces would imply $g$ expands as
    \[g = f^2ds^2 + h.\]
    
    Then, the static equations imply that \[\hess_{S_{\infty}}f = \frac{1}{2}(K + 3)fh,\] where $K$ is the Gaussian curvature of $h$.
    
    In this case, (\cite{Am17}, Proposition 24) implies $h$ has constant Gaussian curvature $ K = -3$. Then, $|\ric(g) + 2g|^2 = 6$ in an open neighborhood of $S_{\infty}$. However, since $|\ric(g) + 2g|^2$ is an analytic function and $g$ is analytic as well, this would imply that $|\ric(g) + 2g|^2 = 6$ everywhere, which contradicts the fact that $\ric(g)\to -2g$ near infinity. Therefore, $ f= 0$ along $S_{\infty}$. 
    
    Lastly, we can cut $(M,g)$ along $S_{\infty}$ to get two static black hole metrics with the same boundary data. So, there is an isometry between the two halves given by the exponential map on the boundary, meaning that by Theorem \ref{thm: bonus}, $(M,g)$ must be a doubling across $S_{\infty}$. 
    \end{proof}

 \section{A family of one ended hyperbolic manifolds that are not optimal for the Penrose Inequality} \label{sect: quasi}
In this section, we give examples of one ended hyperbolic manifolds with outer-minimizing boundary that are not optimal for the Penrose Inequality. The key point is that this family of hyperbolic manifolds is not static. Therefore, they exhibit more flexibility in their scalar curvature which we exploit to show that there are arbitrarily small deformations that increase the area of their outer-minimizing minimal surface boundaries and decrease their mass. 

\subsection{Quasi-Fuchsian Manifolds}\label{sect: sub: quasi}
In this section, we define and establish basic properties of a family of hyperbolic manifolds, called quasi-Fuchsian manifolds. To be precise, 
\begin{definition}(\cite{Sc22}, Definition 1.1)
    A \emph{quasi-Fuchsian} manifold is a complete hyperbolic manifold $M = \R\times \s_{\genus},\genus\geq 2$, and $M$ contains a compact, nonempty, geodesically convex subset. 
\end{definition} 

We explain the natural relation between quasi-Fuchsian manifolds and zero mass Kottler metrics. Recall that the mass zero Kottler metric for a given choice of hyperbolic metric $h_0$ on $\s_{\genus}, \genus\geq 2$ is given by \[g_k = \frac{dr^2}{r^2-1} + r^2h_0.\] In fact, $g_k$ is actually hyperbolic. Using the change of coordinates $r = \cosh(t)$ and doubling $g_k$ over its boundary, we get a complete hyperbolic metric on $\R\times\s_{\genus}$ called a \emph{Fuchsian} that takes the form \[g_f =dt^2 + \cosh^2(t)h_0.\]

Because the Fuchsian metrics are doublings, both of their ends have the same conformal infinity. By Bers Simultaneous Uniformization theorem, all pairs of marked conformal structures at infinity define a unique hyperbolic metric on $\R\times\s_{\genus}$. A \emph{quasi-Fuchsian} manifold is a hyperbolic metric determined by an arbitrary pair of conformal classes on each end of the manifold. In particular, the space of quasi-Fuchsian metrics can be identified with $\R^{6\genus-6}\times\R^{6\genus-6}$ (corresponding to twice the real dimension of Teichm\"uller space).

In Fuchsian metrics, $\{0\}\times\s_{\genus}$ is the unique closed minimal surface. In general, a quasi-Fuchsian manifold can have an arbitrarily large (but finite) number of closed stable minimal surfaces isotopic to $\{0\}\times\s_{\genus}$. A generic property of outer-minimizing surfaces is that they are stable:

\begin{proposition}(Theorem 1.6 and Proposition 8.2, \cite{GuLiVar26})\label{prop: outer_stab}
    An end of a quasi-Fuchsian manifold $(\s_{\genus}\times\R,g)$ admits a minimal surface $\s^\text{out}_{\genus}$ isotopic to $\{0\}\times\s_{\genus}$ which is outermost with respect to that end. Further, there is a total measure set $M^*$ in the space of all quasi-Fuchsian manifolds so that for both ends, $\s^\text{out}_g$ are strictly stable. 
\end{proposition}

\subsubsection{Mass calculation of Quasi-Fuchsian manifolds. }
We describe notions of mass that apply to quasi-Fuchsian manifolds and explain how they fit into the literature. We begin by detailing the system of coordinates that we will use for the remainder of the section. 

Let $g_f$ be a Fuchsian metric. Let $h_0$ be the hyperbolic metric in the conformal class at infinity of $g_f$. By applying a change of coordinates $x = 2e^{-t}$ to $g_f$,  
    \[g_f = \frac{1}{x^2} \left(dx^2+\left(1+\frac{x^2}{4}\right)^2h_0\right).\] By comparing to the Fefferman—Graham expansion of $g_f$: \[g_f = \frac{dx^2 + h_0+x^2h_2^f + x^4h_4^f}{x^2},\] we see that $h_2^f$ and $h_4^f$ are tensors on $\s_{\genus}$ and $h_2^f = h_0/2,h_4^f = h_0/16$.
    
    Let $g_q$ be a quasi-Fuchsian metric. If the conformal infinity at an end of $g_q$ is $[h_0]$, then near that end, the Fefferman—Graham expansion of $g_q$ is given by \[g_q = \frac{dx^2 + h_0+x^2h_2^q + x^4h_4^q}{x^2},\] where $h_2^q$ and $h_4^q$ are again tensors on $\s_{\genus}$. 

    Further, standard properties of expansions of Einstein metrics imply \[\tr_{h_0}(h_2^f-h_2^q) = 0\] because they have the same conformal infinity and that $\tr_{h_0}h_2 = -\frac{1}{2}R(h_0)$. And so also, \[\Div_{h_0} (h_2^q-h^f_2) = d(\tr_{h_0} (h_2^q-h_2^f)) = 0.\] 

First, we recall the Chru\'sciel-Herzlich definition of mass in the case where the reference manifold $(\model,\mg)$ is Fuchsian: $\model = \R\times\s_{\genus}$ where $\genus\geq 2$. Note that the family of static potentials on $\model$ is one dimensional. 

\begin{definition}\label{def: mass}
Let $g$ be a manifold ALH to a Fuchsian metric $g_f$ at rate a rate $q > 3/2$ with chart at infinity $\Phi$. Let $e = \Phi^*g - g_f$. Let $\nu = -\frac{1}{x}\del_x$ be the unit normal to the surface $x = s$. Let norms, $\grad$, $\Div$, and $\tr$ be taken with respect to $g_f$. Let $V(x) = \frac{1}{x} - \frac{x}{4}$ the unique static potential (up to scaling) for $g_f$. 
\[
\int_{M_{ext}}
\left( |e|^2 + \left|\nabla e\right|^2 \right)
V(x) \, dg_f
< \infty,
\]
and
\[
\int_{M_{ext}}
\left|R_g + n(n-1)\right|
V(x)\, dg_f
< \infty,
\]
then the following integral exists and is finite:
\[16\pi(\genus-1)m(g) = \lim_{s\to 0}\int_{x = s} (V(\Div e-d(\tr e)) - e(\grad V,\cdot)+\tr (e)dV)(\nu)dS.\] We call $m(g)$ the \emph{mass} of $(M,g)$. 
\end{definition}
Note that we chose the normalization by $16\pi(\genus-1)$ so that a Kottler metric with parameter $m$ has mass $m$. 
 
\begin{proposition}\label{prop: quasi has 0 mass}
    The mass of a quasi-Fuchsian manifold $(M, g_q)$ is 0. 
\end{proposition}
\begin{proof}
   For a quasi-Fuchsian metric $g_q$, the reference metric chosen is the Fuchsian metric $g_f$ with the same conformal infinity at that end. This way $|g_q-g_f| = O(x^2)$ and the mass is well defined. 

    Let $ e= g_q-g_f.$ Then, $e_{xx} = 0$ and $e_{xi} = 0$. Note that $\nu = -x\del_x$. 
    Note that \[g_f^{ij} = O(x^2), e_{aa} = O(1), \text{ and }\Gamma_{ax}^b = (-\frac{1}{x} + O(x))\delta_a^b.\] 

    Since $\tr_{h_0}(h_2^q-h^f_2) = 0$ and $\tr_{h_0} e = O(x^2)$,  \[\tr e = x^2\left(1+\frac{x^2}{4}\right)^{-2}\tr_{h_0}e =O(x^2)\cdot O(x^2)=O(x^{4}).\] Since additionally $(d\tr_{h_0}e)(\nu) =-x\del_x\tr_{h_0} e= O(x^2)$,  \[d\tr (e)(\nu) = -x\del_x\big(x^2(1+\frac{x^2}{4})^{-2}\tr_{h_0}(e)\big) = O(x^4).\] 
    Next, using that $e_{xa} = 0$, $e_{xx} = 0$, and $\tr_{g_f} e= O(x^4)$,
    \begin{align*}
        \Div (e)(\nu) &= -xg_f^{ij}(\del_je_{xi} -\Gamma_{jx}^me_{mi} - \Gamma_{ji}^me_{xm})\\
        &= xg_f^{ij}\Gamma_{jx}^me_{mi}\\
        &= xg_f^{ij}\Gamma_{jx}^je_{ji}\\
        &= (-1+O(x^2))\tr_{g_f}e\\
        &= O(x^4).
    \end{align*}
    
    Since $V= (\frac{1}{x}-\frac{x}{4}) = O(x^{-1})$ and $dV(\nu) = \frac{1}{x} +\frac{x}{4} = O(x^{-1})$, so \[\tr (e)dV(\nu) = O(x^3)\] and \[V(\Div (e)-d\tr(e))(\nu) = O(x^3).\]

    Lastly, $e(\grad V,\nu) = 0$ because $e_{\nu\nu}  =0.$ Therefore, for some sufficiently small $x_0$, the mass integral can be written as \[\int_{x = s} (V(\Div e-d(\tr e)) - e(\grad V,\cdot)+\tr (e)dV)(\nu)dS =\int_{x = x_0} O(s)dS\] which goes to 0 as we take $s\to 0$, as desired. 
\end{proof}

\begin{lemma}\label{lem: mass is neg II}
Let $u = 1+u_3x^3 + f$ where $f = O(x^4)$. If $u_3\leq 0$ and $u_3$ is not identically zero, the metric $g = u^4{g_q}$ has negative mass.
\end{lemma}
\begin{proof}
    We proceed similarly to the proof of Proposition \ref{prop: quasi has 0 mass}. Let $e = g - g_f.$ Note that \[e = g - g_f = (1+u_3x^3 + f)^4g_q-g_f =g_q-g_f +4u_3x^3g_f + O(x^4).\]  
    Let $\varphi = 4u_3x^3$. Note that \[\Div(\varphi g_f) = d\varphi\] and \[\tr(\varphi g_f) = 3\varphi.\] Since $m(g_q) = 0$ and any terms $O(x^4)$ or higher do not contribute to the mass, the mass integral can be computed by: \[16\pi(\genus-1)m(g) = \lim_{s\to 0}\int_{x = s} -2Vd\varphi(\nu) + 2\varphi dV(\nu) = \lim_{s\to 0}\int_{x = s}32u_3x^2dS < 0.\] 
\end{proof}

\subsection{Non-Fuchsian quasi-Fuchsian metrics are not static}

In this subsection, we prove that quasi-Fuchsian manifolds are not static. We also show that neighborhoods of separating sets of quasi-Fuchsian manifolds are not static. 

\begin{proposition}\label{prop: quasi-Fuchsian not static}
Non-Fuchsian quasi-Fuchsian manifolds are not static.
\end{proposition}
\begin{proof}  
    According to Definition \ref{def:alh}, either end of a quasi-Fuchsian metric $g$ is ALH at rate $q = 2$. Therefore, by Theorem \ref{thm:main}, if a quasi-Fuchsian manifold were static, it would be the doubling of an ALH static black hole metric. In particular, the doubled horizon would give a compact totally geodesic hypersurface for $g$. By (Theorem 3.3, \cite{Uh83}), the existence of a closed totally geodesic hypersurface implies $g$ is Fuchsian, contrary to assumption. 
\end{proof}

We want to be able to say that particular subsets of a quasi-Fuchsian manifold $M$ are also not static. It is not true that \emph{any} subset of $M$ is static, even when $M$ is not static. For example, a simply connected ball in $M$ has $4$ static potentials because it can be embedded in $\mathbb{H}^3$. Another potential issue is that some $S\subset M$ could support a static potential that cannot be extended to the entire manifold. 

\nopmt*
\begin{proof}
     By (\cite{corvino2000}, Theorem 4), if $\Omega$ is an open set smooth, relatively compact such that $(\Omega, g)$ is not static, then its scalar curvature can be prescribed, granted the prescription is close enough to $R(g)$. 
     
     Suppose that $(\Omega,g)$ supports a non-trivial static potential $f$. Consider the lift $\tilde{\Omega}$ of $\Omega$ and $\tilde{f}$ of $f$ to $\mathbb{H}^3$. Since the space of solutions to the static equations has maximal dimension, $\tilde{f}$ must extend globally to a static potential $\bar{f}$ on $\mathbb{H}^3$. Because $\pi_1(\Omega)\twoheadrightarrow \pi_1(M)$, $\bar{f}$ must descend to all of $M$ so the potential $f$ originally supported on $\Omega$ is actually defined on all of $M$. However, by Proposition \ref{prop: quasi-Fuchsian not static}, a quasi-Fuchsian manifold which is not Fuchsian is not static, so this is a contradiction.
\end{proof}

\subsection{Construction of the Examples.}

\begin{lemma}\label{lem: mass is neg}
    Let $g_{\e}$ be a metric on $\R\times\s_{\genus}$ that is equal to a quasi-Fuchsian metric $g$ (which is not Fuchsian) outside a compact set, $R(g_{\e})\geq -6$ and $R(g_{\e})>-6$ somewhere. Then, near either end, the solution to the Yamabe problem $u$ expands as $u = 1 + u_3x^3 + O(x^4)$, $u_3 \leq 0$, and $u_3$ is not identically zero.  
\end{lemma}
\begin{proof}    
    Let $\del_{\infty} M$ be the conformal infinity of an end of $M$ and let $x$ be a defining function near $\del_{\infty} M$. From (Theorem 1.3 \cite{AnChFr92}), a solution $u$ to the Yamabe problem expands asymptotically as \[ u = 1 + \sum_{i = 1}^{\infty}\sum_{j = 0}^{N_i}u_{ij}(\log x) ^jx^i,\] where $u_{ij}\in C^{\infty}(\del_{\infty} M)$ and $N_i\in\N$ and $N_i = 0$ when $i<3$. Because $g_{\e}$ is quasi-Fuchisan outside a compact set, using the analysis from (Theorem 1.3.2 and Lemma 2.1 \cite{AnChFr92}) implies that $u$ actually expands as \[ u = 1 + \sum_{i = 1}^{\infty}u_{ij}(\log x) ^jx^i,\] where $u_{ij}\in C^{\infty}(\del_{\infty} M)$.

    We show $u_3\leq 0$ and is not identically zero. Let $u = 1+ u_mx^m + O(x^{m+1})$, where $m\geq 3$ is the first $m$ such that $u_m\not\equiv 0$. To start, $R(g_{\e}) \geq -6$ implies $u < 1$ by the maximum principle, so $u_m\leq 0$. Using the Yamabe equation:
    \begin{align*}
        -8\lap_g u -6u &= -6 u^5\\
        -8\lap(1+u_mx^m) - 6(1+u_mx^m) &= -6(1+5u_mx^m) + O(x^{m+1})\\
        \lap_gu_mx^m - 3u_mx^m &= O(x^{m+1}).
    \end{align*}
    Since $\lap (u_mx^m) = u_m\lap x^m + x^m\lap u_m, \lap x^m = m(m-2)x^m + O(x^{m+1})$, and $\lap u_m = O(x^2)$, \[\lap_g (u_mx^m) - 3u_mx^m = (m(m-2) - 3)u_mx^m + O(x^{m+2})\]if $m > 3$, $\lap_gu_mx^m - 3u_mx^m = O(x^{m})$. Which implies that if $u_3\equiv 0$, then $u-1$ vanishes to infinite order at infinity.  However, $u-1$ cannot vanish to infinite order at the boundary: $u-1$ satisfies the equation \[4\lap_g (u-1) -3(u+1)(u^2+1)(u-1) = 0\] far enough into either end. But, then (\cite{Ma91}, (14) Theorem) would imply that $u-1\equiv 0$, which is a contradiction. Therefore, $u_3$ is not identically 0, as desired.     
\end{proof}

Fix a quasi-Fuchsian manifold $(M_q,g_q)$ and a particular end of $(M_q,g_q)$. Assume that the outermost minimal surface $\s^{\text{out}}$ with respect to that end is strictly stable, which exists on a full measure set of quasi-Fuchsians by \ref{prop: outer_stab}. Let $(M,g)$ be an end of $(M_q,g_q)$ cut along $\s^{\text{out}}$ containing no other closed minimal surfaces. 

Let $A$ be the area of $\s^{\text{out}}$. Recall that the mass of $(M,g)$ is 0 by Proposition \ref{prop: quasi has 0 mass}. We recall the main theorem of this section:

\nopenroseqf*
\begin{proof}
    Let $\s$ be a strictly stable outer-minimizing surface in $M$. Let $S$ be a smooth separating surface in $M$ such that $\dist(\s,S)> 1$. Let $S_{t}$ be the equidistant surface of distance $t\in\R$ from $S$. Assume that the injectivity radius of $S$ in $M$ is also greater than one so that $S_1$ and $S_{-1}$ are non-singular and separating. Cutting $M$ along $S_{-1}$ and $S_1$ splits $M$ into three pieces: $M_1$, which contains $\s$, $M_2$, $M_2$ which is bounded by $S_{-1}$ and $S_1$, and $M_3$ which is the unbounded end containing no closed minimal surfaces. See the figure below for a diagram. \\
    
    \noindent {\bf Increasing the area of the outer-minimizing minimal surface.} Let $\di > 0$ and let $\eta_{\di}\colon M\to \R_{\geq 0}$ be a bump function compactly supported in $M_1$ such that the support of $\eta_{\di}$ does not intersect $\s$ and such that $|\eta_{\di}|_{C^{4,\alpha}}\leq \di$. By Theorem \ref{lem: quasi-Fuchsian scal curv}, there is a perturbation $g_{\di}^{\p}$ of $g$, compactly contained in $M_1$ such that $R(g_{\di}^{\p}) = -6-\eta_{\di}$. Since $\eta_{\di}$ is not supported on $\s$, we may assume that $g = g_{\di}^{\p}$ in a small neighborhood of $\s$.  

    Since $\del M$ was a stable nondegenerate minimal surface with respect to $g$, it is a local area minimizer with respect to $g$. 
    Since $\eta_{\di}$ is supported away from $\del M$, we can shrink $\di$ if necessary so that $\s$ remains outer-minimizing minimal for $g_{\di}^{\p}$.
    
    By (\cite{AnChFr92} Theorem 3.4), we can apply a conformal change of metric by the function $u_{\di}^4$ to $g_{\di}^{\p}$ to get a metric $g_{\di}\equiv u_{\di}^4g_{\di}^{\p}$ such that \[R(g_{\di}) = -6 \text{ and } u_{\delta_1}\to 1\text{ as } |x|\to\infty.\]  Recalling that $\s$ is strictly stable and using an implicit function type argument as in \cite{huang22}, we see that as $|g_{\di}^{\p}-g|\to 0$ in $C^{2,\alpha}_{-q}$, $|u_{\delta_1}-1|\to 0$ in $C^{2,\alpha}_{-q}$. Therefore, by shrinking $\delta_1$ if needed
 we may assume that $g_{\di}$ still has a closed outer-minimizing minimal surface $\s_{\di}$. By the maximum principle, $u_{\di} >1$. Therefore, \[\text{Area}(\s_{\delta_1},g_{\di}) > \text{Area}(\s,g_{\di}^{\p}) = \text{Area}(\s,g).\] \\

 \noindent {\bf Decreasing the mass.} Now let $\dii > 0$ and let $\eta_{\dii}\colon M\to\R_{\geq 0}$ be a smooth bump function compactly contained in $M_3$ so that $|\eta_{\dii}|_{C^{4,\alpha}}\leq \dii$. Again by Theorem \ref{lem: quasi-Fuchsian scal curv}, there is a perturbation $g_{\dii}^{\p}$ of $M$ compactly contained on $M_3$ such that $R(g_{\dii}^{\p}) = -6+\eta_{\dii}$.
    Again, by (\cite{AnChFr92} Theorem 3.4), we can apply a conformal change of metric by the function $u_{\dii}^4$ to $g_{\dii}^{\p}$ to get a metric $g_{\dii}\equiv u_{\dii}^4g_{\dii}^{\p}$ such that \[R(g_{\dii}) = -6 \text{ and } u_{\delta_2}\to 1\text{ as } |x|\to\infty.\]   By Lemma \ref{lem: mass is neg} and Lemma \ref{lem: mass is neg II}, the mass of $(M,g_{\dii})$ is negative. Using again the same implicit function type argument as in \cite{huang22}, we see that as $|g_{\dii}^{\p}-g|\to 0$ in $C^{2,\alpha}_{-q}$, $|u_{\delta_2}-1|\to 0$ in $C^{2,\alpha}_{-q}$. Therefore, by shrinking $\delta_2$ if needed we may assume that $g_{\dii}$ has no closed outer-minimizing minimal surfaces intersecting $M_{2}\cup M_3$. \\    

    \noindent {\bf Gluing.} Lastly, we glue $(M_1\cup M_2,g_{\di})$ and $(M_2\cup M_3,g_{\dii})$ together across $M_2$. Consider a smooth function $\Phi\colon M\to \R_{\geq 0}$ such that $\Phi(x) =0$ in a neighborhood of $M_1$ and $\Phi(x) =1$ in a neighborhood of $M_3$. Then we can define the metric 
    \[g^{\p} = 
    \begin{cases}
        g_{\di} & x\in M_1\\
        (1-\Phi(x))g_{\di} + \Phi(x)g_{\dii} & x\in M_2\\
        g_{\dii} & x\in M_3.
    \end{cases}
    \]
 Since the estimates from Theorem 1 in \cite{corvino2000} can be taken uniform near $g$, shrinking $\di$ and $\dii$ sufficiently small gives a compact perturbation $\bar{g}$ of $g^{\p}$ in $M_2$ so that $R(\bar{g}) = -6$ such that $\s_{\di}$ is the outer-minimizing minimal surface. 
    
    We can cut $(M,\bar{g})$ along $\s_{\di}$ to obtain $(\bar{M},\bar{g})$ as the unbounded component. Shrinking $\di,\dii\to 0$ shows that $|g-\bar{g}|_{C^{2,\alpha}}$ can be taken arbitrarily small. 

    \begin{tikzpicture}
      \hspace{.125\linewidth}  \node (img) {\includegraphics[width=.75\linewidth]{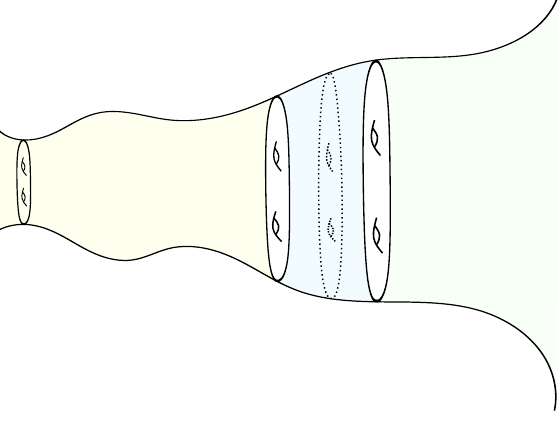}};
\node (S) at (1.1,3.25) {$S$};
\node (plus) at (2.1,3.55) {$S_1$};
\node (minus) at (-0.15,2.8) {$S_{-1}$};
\node (minsur) at (-5.5,1.85) {$\s$};
\draw [decorate,decoration={brace,amplitude=5pt,mirror,raise=4ex}]
  (-6,-3.65) -- (-0.02,-3.65) node[midway,yshift=-3em]{$M_1$};
\draw [decorate,decoration={brace,amplitude=5pt,mirror,raise=4ex}]
  (0.01,-3.65) -- (2.2,-3.65) node[midway,yshift=-3em]{$M_{2}$};
  \draw [decorate,decoration={brace,amplitude=5pt,mirror,raise=4ex}]
  (2.23,-3.65) -- (6,-3.65) node[midway,yshift=-3em]{$M_3$};
\end{tikzpicture}
\end{proof}

\bibliographystyle{abbrv}
{\footnotesize\bibliography{bib}}

\end{document}